\newif\ifpagetitre            \pagetitretrue
\newtoks\hautpagetitre        \hautpagetitre={\hfil}
\newtoks\baspagetitre         \baspagetitre={\hfil}
\newtoks\auteurcourant        \auteurcourant={\hfil}
\newtoks\titrecourant         \titrecourant={\hfil}

\newtoks\hautpagegauche       \newtoks\hautpagedroite
\hautpagegauche={\hfil\the\auteurcourant\hfil}
\hautpagedroite={\hfil\the\titrecourant\hfil}

\newtoks\baspagegauche \baspagegauche={\hfil\tenrm\folio\hfil}
\newtoks\baspagedroite \baspagedroite={\hfil\tenrm\folio\hfil}

\headline={\ifpagetitre\the\hautpagetitre
\else\ifodd\pageno\the\hautpagedroite
\else\the\hautpagegauche\fi\fi}

\footline={\ifpagetitre\the\baspagetitre
\global\pagetitrefalse
\else\ifodd\pageno\the\baspagedroite
\else\the\baspagegauche\fi\fi}

\vsize=9.0in\voffset=1cm
\looseness=2


\message{fonts,}

\font\tenrm=cmr10
\font\ninerm=cmr9
\font\eightrm=cmr8
\font\teni=cmmi10
\font\ninei=cmmi9
\font\eighti=cmmi8
\font\ninesy=cmsy9
\font\tensy=cmsy10
\font\eightsy=cmsy8
\font\tenbf=cmbx10
\font\ninebf=cmbx9
\font\tentt=cmtt10
\font\ninett=cmtt9

\font\ninesl=cmsl9
\font\eightsl=cmsl8

\font\nineit=cmti9
\font\eightit=cmti8

\skewchar\ninei='177 \skewchar\eighti='177
\skewchar\ninesy='60 \skewchar\eightsy='60

\def\eightpoint{\def\rm{\fam0\eightrm} 
\normalbaselineskip=9pt
\normallineskiplimit=-1pt
\normallineskip=0pt

\textfont0=\eightrm \scriptfont0=\sevenrm \scriptscriptfont0=\fiverm
\textfont1=\ninei \scriptfont1=\seveni \scriptscriptfont1=\fivei
\textfont2=\ninesy \scriptfont2=\sevensy \scriptscriptfont2=\fivesy
\textfont3=\tenex \scriptfont3=\tenex \scriptscriptfont3=\tenex
\textfont\itfam=\eightit  \def\it{\fam\itfam\eightit} 
\textfont\slfam=\eightsl \def\sl{\fam\slfam\eightsl} 

\setbox\strutbox=\hbox{\vrule height6pt depth2pt width0pt}%
\normalbaselines \rm}

\def\ninepoint{\def\rm{\fam0\ninerm} 
\textfont0=\ninerm \scriptfont0=\sevenrm \scriptscriptfont0=\fiverm
\textfont1=\ninei \scriptfont1=\seveni \scriptscriptfont1=\fivei
\textfont2=\ninesy \scriptfont2=\sevensy \scriptscriptfont2=\fivesy
\textfont3=\tenex \scriptfont3=\tenex \scriptscriptfont3=\tenex
\textfont\itfam=\nineit  \def\it{\fam\itfam\nineit} 
\textfont\slfam=\ninesl \def\sl{\fam\slfam\ninesl} 
\textfont\bffam=\ninebf \scriptfont\bffam=\sevenbf
\scriptscriptfont\bffam=\fivebf \def\bf{\fam\bffam\ninebf} 
\textfont\ttfam=\ninett \def\tt{\fam\ttfam\ninett} 

\normalbaselineskip=11pt
\setbox\strutbox=\hbox{\vrule height8pt depth3pt width0pt}%
\let \smc=\sevenrm \let\big=\ninebig \normalbaselines
\parindent=1em
\rm}

\def\tenpoint{\def\rm{\fam0\tenrm} 
\textfont0=\tenrm \scriptfont0=\ninerm \scriptscriptfont0=\fiverm
\textfont1=\teni \scriptfont1=\seveni \scriptscriptfont1=\fivei
\textfont2=\tensy \scriptfont2=\sevensy \scriptscriptfont2=\fivesy
\textfont3=\tenex \scriptfont3=\tenex \scriptscriptfont3=\tenex
\textfont\itfam=\nineit  \def\it{\fam\itfam\nineit} 
\textfont\slfam=\ninesl \def\sl{\fam\slfam\ninesl} 
\textfont\bffam=\ninebf \scriptfont\bffam=\sevenbf
\scriptscriptfont\bffam=\fivebf \def\bf{\fam\bffam\tenbf} 
\textfont\ttfam=\tentt \def\tt{\fam\ttfam\tentt} 

\normalbaselineskip=11pt
\setbox\strutbox=\hbox{\vrule height8pt depth3pt width0pt}%
\let \smc=\sevenrm \let\big=\ninebig \normalbaselines
\parindent=1em
\rm}

\message{fin format jgr}

\hautpagegauche={\hfill\ninerm\the\auteurcourant}
\hautpagedroite={\ninerm\the\titrecourant\hfill}
\auteurcourant={R.G.\ Novikov}
\titrecourant={A holographic uniqueness theorem}

\magnification=1200
\font\Bbb=msbm10
\def\R{\hbox{\Bbb R}}

\def\S{\hbox{\Bbb S}}
\def\pa{\partial}

\vskip 2 mm
\centerline{\bf A holographic uniqueness theorem}

\vskip 2 mm
\centerline{\bf R.G.\ Novikov}

\vskip 4 mm
\noindent
{\bf Abstract.}
We consider a plane wave, a radiation solution, and the sum of these solutions (total solution)
for the Helmholtz equation in an exterior region in $\R^3$. 
We consider a ray in this region, such that its direction is different from the propagation direction of the plane wave.
We show that the restriction of the radiation solution  to this ray is uniquely determined
by the intensity of  the  total solution on an interval of this ray.
As a corollary, we also obtain that  the restriction of the radiation solution to any plane in the exterior region is uniquely determined  by  the  intensity 
of  the  total solution on an open domain in this plane.
In particular, these results solve one of old mathematical questions of holography.


\vskip 2 mm
{\bf 1. Introduction}

We consider  the Helmholtz equation
$$ -\Delta\psi(x)=\kappa^2 \psi(x),\ \ x\in {\cal U},\ \  \kappa>0, \eqno(1.1)$$
where $\Delta$ is the Laplacian in $x$, and $\cal U$ is a region (open connected set) in $\R^3$
consisting of all points outside a closed bounded regular surface $S$ (as in [A], [W]).

This equation arises, in particular, in electrodynamics, acoustics, and quantum mechanics.

For equation (1.1) we consider solutions  $\psi_0$ and $\psi_1$ such that:

$$\psi_0= e^{ikx},\ \ k\in\R^3,\ \ |k|=\kappa ;   \eqno(1.2)$$
$\psi_1$ is of class $C^2$ and satisfies the Sommerfeld's radiation condition
$$  |x|({\pa\over \pa|x|}-i\kappa)\psi_1(x)\to 0\ \ as \ |x|\to +\infty, \eqno(1.3)$$
uniformly in $x/|x|$. 
We say that $\psi_0$  is a plane wave solution and $\psi_1$  is a radiation solution. 

Let  $L=L_{x_0, \theta}$ be the ray in  $\R^3$  that starts at $x_0$ and has direction  $\theta$:

$$ L=L_{x_0, \theta}=\{x\in\R^3:\ \ x=x(s)=x_0+s\theta,\ \ 0<s <+\infty  \}, \ \ x_0\in\R^3, \ \  \theta\in\S^2,  \eqno(1.4) $$
where $\S^2$ is the unit sphere in $\R^3$.

In the present work we show that, for a fixed  plane wave solution $\psi_ 0$, 
any complex-valued  radiation solution $\psi_ 1$  on $L$ is uniquely determined  by  the intensity $|\psi|^2$  of the total solution $\psi= \psi_ 0+\psi_ 1$ on an arbitrary interval of $L$,
under the assumptions that $L=L_{x_0, \theta}\subset ~{\cal U}$ and $\theta \neq k/|k|$. Here,  ${\cal U}$ is the region in (1.1) and $k$ is the vector in (1.2).  
This result is given as Theorem 2.1 in Section 2.

As a corollary, we also obtain that, for any plane $X\subset  {\cal U}$, any complex-valued  $\psi_ 1$  on $X$ is uniquely determined  by  the  intensity $|\psi_ 0+\psi_ 1|^2$ 
on an arbitrary open domain of $X$,  for  fixed   $\psi_ 0$; see Theorem 2.2 in Section 2.

Our studies are motivated by holography  and phaseless inverse scattering; see, e.g.,  [G], [Wo1],  [Wo2], [JL], [ Nu],  [Kl1], [M], [N1]-[N3], [KR], [MH], [Kl2], [R], [HN], [NS2]
and references therein.
As far as we know,  in the literature the aforementioned determination of $\psi_1$ on  $X$ from  $|\psi_ 0+\psi_ 1|^2$  on  $X$ was considered in the framework of different approximations only.
Therefore, Theorems 2.1 and 2.2 solve one of old mathematical  questions arising in holography.
In addition, various new results on phaseless inverse scattering follow from Theorems 2.1 and 2.2 and known results on inverse scattering with phase information.

The main results of this work are presented in more detail and proved in Sections 2,  4, and 5.
In our proofs we proceed from results recalled in Section 3.

\vskip 2 mm
{\bf 2. Main results}

Our key result is as follows.

\vskip 2 mm
{\bf Theorem 2.1.} {\sl Let $\psi_0$ and  $\psi_1$  be solutions of equation (1.1) as in formulas (1.2) and (1.3).
Let $L$ be a ray as in (1.4) such that  $L=L_{x_0, \theta}\subset {\cal U}$,  $\theta \neq k/|k|$, where ${\cal U}$ is the region in (1.1), $k$ is the wave vector in (1.2).
Then $\psi_1$ on $L$ is  uniquely determined by the  intensity  $|\psi_0+\psi_1|^2$ on   $\Lambda$,  for fixed $k$,
where   $\Lambda$  is an arbitrary non-empty open interval of~$L$.}

\vskip 2 mm
As a corollary, we also get the following result.

\vskip 2 mm
{\bf Theorem 2.2.} {\sl Let $\psi_0$ and  $\psi_1$  be solutions of equation (1.1) as in formulas (1.2) and (1.3).
Let $X$ be a two-dimensional plane in $\R^3$ such that $X\subset {\cal U}$.
Then $\psi_1$ on $X$ is  uniquely determined by the  intensity  $|\psi_0+\psi_1|^2$ on  ${\cal D}$,  for fixed $k$,
where  ${\cal D}$  is an arbitrary non-empty open domain of $X$.}

\vskip 2 mm
Theorem 2.1 is proved in Sections 4 and 5 using techniques developed in [A], [W],  [N1], [N2], [N4], [NS2],  [NSh].

Note that $\psi_0$, $\psi_1$ are real-analytic on ${\cal U}$, and, therefore, on $L$ in Theorem 2.1 and on $X$ in Theorem 2.2.
Therefore,  the function $|\psi|^2=|\psi_0+\psi_1|^2= (\psi_0+\psi_1)(\overline \psi_0+\overline \psi_1)$ is  real-analytic on $L$ in Theorem 2.1 and on $X$ in Theorem 2.2.
In view of this analyticity,  Theorem~2.1 reduces to the case when  $\Lambda=L$ and Theorem 2.2 reduces to the case when ${\cal D}=X$.

Theorem 2.2 is proved as follows (for example).
We assume that ${\cal D}=X$.  
Then to determine $\psi_1$ at  an arbitrary $x\in X$,  we consider a ray $L=L_{x_0, \theta}\subset X$  such that $x\in L$,  $\theta \neq k/|k|$,  and  use Theorem 2.1 for this $L$.

\vskip 2 mm
{\bf Remark 2.1.} In the one dimensional case, an analog of Theorem 2.1 follows from results of [N3].
The two dimensional case will be considered elsewhere using  results of [K] in place of the three dimensional Atkinson-Wilcox expansion (3.2) used in the present work.

To our knowledge, in the literature,  rigorous mathematical results on the determination of the radiation solution $\psi_ 1$ on  $X$ from the intensity   $|\psi|^2$ 
of the total solution  $\psi=\psi_ 0+\psi_ 1$  on  $X$  (including non-uniqueness and uniqueness),
are given only in the framework of  considerable approximations (e.g., paraxial approximation  for the case when $k$ is orthogonal to $X$); see [Nu], [M], [MH].
Therefore,  Theorem 2.2  solves  an old mathematical question of holography.

In fact,  an asymptotic  determination of  $\psi_ 1$ on $L_{x_0, \theta}$  from the intensity   $|\psi_ 0+\psi_ 1|^2$ on $L_{x_0, \theta}$ was developed in [N1],  [N2],  [N4], [NS1], [NS2],  [NSh],
for $x_0=0$.  In Section~3 we recall formulas (3.5)-(3.8)  in this connection.
In the present work we proceed from these asymptotic results and classical results on the radiation solution $\psi_ 1$, including the Atkinson-Wilcox expansion (3.2).

Many new uniqueness results on phaseless inverse (coefficient and obstacle) scattering problems follow from Theorems 2.1 and 2.2 and known results on inverse scattering with phase information
in three dimensions. These new uniqueness results will be listed  elsewhere.

Note  that the determination in Theorems 2.1 and 2.2 includes an analytic continuation. 
Studies on more stable reconstructions appropriate for numerical  implementation  will be continued  elsewhere.

\vskip 2 mm
{\bf 3. Preliminaries}

Let

$$B_r=\{x\in\R^3:\ \ |x|<r\},\ \ r>0.\eqno(3.1)$$

Suppose that  $\psi_1$ is a radiation solution of (1.1), and  $\R^3\setminus B_r  \subset {\cal U}$. Then, due to results of [A], [W],  we have that 

$$\psi_1(x)={e^{i\kappa|x|}\over |x|}\sum\limits_{j=1}^{\infty}{f_j(\theta)\over |x|^{j-1}} \ \ {\rm for}\ \ x\in \R^3\setminus B_r,\ \ \theta={x\over |x|},\eqno(3.2)$$

where the series converges absolutely and uniformly.

Let 

$$a(x,k)=|x|(|\psi(x)|^2-1),\ \ x\in{\cal U}, \eqno(3.3)$$

where $\psi=\psi_0+\psi_1$,  $\psi_0$  and $\psi_1$ are  solutions of (1.1) as in (1.2) and (1.3), $k$ is the wave vector in (1.2).
Then (see [N1], [N4], [NSh]):

$$a(x,k)=e^{i(\kappa|x|-kx)}f_1(\theta)+  e^{-i(\kappa|x|-kx)} \overline {f_1(\theta)}+O(|x|^{-1})\ \ as \ |x|\to +\infty,  \eqno(3.4)$$

uniformly in $\theta=x/|x|$;

$$f_1(\theta)={1\over D}\bigl(e^{i(ky-\kappa|y|)}a(x,k)-e^{i(kx-\kappa|x|)}a(y,k)+O(|x|^{-1}) \bigr)\ \ as \ |x|\to +\infty,  \eqno(3.5)$$
$$D=2isin(\tau(k\theta-\kappa)),\ \  \theta\in\S^2,\ \  \tau>0, $$
where
$$ x,y\in L_{x_0, \theta},\ \ x_0=0,\ \ y=x+\tau\theta,  \eqno(3.6)$$
and $D\ne 0$ for fixed $\theta$ and  $\tau$.

\vskip 2 mm
{\bf Remark 3.1.} If an arbitrary function $a$ on  $L_{0, \theta}$ satisfies (3.4),  for fixed $\theta\in\S^2$,  $k\in\R^3$,  $\kappa=|k|>0$,
then formula (3.5) holds.

\vskip 2 mm
Moreover, in fact, our work [N4] gives  (see formulas (3.4) in  [N4]):
$$\eqalignno{
&{\rm formulas\ for\ finding}\ \ f_j(\theta)\ \ {\rm up\ to}\ \ O(s^{-(n-j+1)})\ \ {\rm as}\ \  s\to +\infty,\ j=1,\ldots,n,&(3.7)\cr
&{\rm from}\ \ |\psi(x)|^2\ \ {\rm given\ at}\ \  2n\ \ {\rm points}\ \ x=x_1(s),\ldots,x_{2n}(s)\ \ {\rm on}\ \   L_{0, \theta},\cr}$$
where
$$\eqalign{
&x_i(s)=r_i(s)(\theta),\ \ i=1,\ldots,2n,\cr
&r_{2j-1}(s)=\lambda_js,\ \ r_{2j}(s)=\lambda_js+\tau,\ \ j=1,\ldots,n,\cr
&\lambda_1=1,\ \ \lambda_{j_1}<\lambda_{j_2}\ \ {\rm for}\ \ j_1<j_2,\ \
 \tau>0,\cr}\eqno(3.8)$$
and $D\ne 0$ for fixed $\theta$ and  $\tau$, where $\psi$ is as in (3.3),  $D$ is as in (3.5).

For $n=1$, formulas (3.7), (3.8) reduce to (3.5). In general, formulas (3.7), (3.8) are given recurrently in [N4].
Under the assumption that $\theta \neq k/|k|$, the parameter $\tau$ can be always fixed in such a way that  $D\ne 0$.

\vskip 2 mm
{\bf 4. Proof of  Theorem 2.1}

\vskip 2 mm
{\it 4.1. Case  $L\subseteq L_{0, \theta}$.} First, we give the proof for the case when $L\subseteq L_{0, \theta}$.
In this case it is sufficient that  $|\psi|^2$ on $L$ uniquely determines  $f_j(\theta)$  in (3.2) for all $j$ via  formulas (3.7), (3.8),  where $\psi=\psi_0+\psi_1$.
The rest follows from the convergence of the series in (3.2) and analyticity of $\psi_1$ and  $|\psi|^2$  on~$L$.

Note that the recurrent formulas (3.7), (3.8) are not very simple in [N4]. 
Therefore, for completeness of presentation in Appendix,  we  give a very simple proof that $|\psi|^2$ on $L$ uniquely determines  $f_j(\theta)$  in (3.2) for all $j$.

\vskip 2 mm
{\it 4.2. General case.} In fact, the general case reduces to the case of Subsection 4.1 by the change of variables 
$$x'=x-q,\ \ {\rm  such\ that}\ \  L\subseteq L_{q, \theta},  \eqno(4.1)$$
for some fixed $q\in\R^3$.

In the new variables $x'$, we have that:
$$\psi_0= e^{ikq}e^{ikx'}, \eqno(4.2)$$  
$$\eqalign{
&\psi_1\ \ {\rm   satisfies}\ \  (3)\ \ {\rm and\ admits\ presentation}\ \  (3.2)\ \  {\rm with}\ \ x'\ \  {\rm in\ place\ of}\ \ x, \cr
&{\rm with\ some\ new}\ \  f_j'\ \ {\rm and }\ \ r'\ \  {\rm in}\ \ (3.2),     \cr}  \eqno(4.3) $$  
$$|\psi|^2=|e^{ikx'}+ e^{-ikq}\psi_1(x')|^2, \eqno(4.4)$$
where $x'\in{\cal U'}={\cal U}-q$,  and $r'$ is such that  $\R^3\setminus B_{r'}  \subset {\cal U'}$;
$$L\subseteq L_{q, \theta}=L_{0, \theta}.  \eqno(4.5)$$.

In addition,
$$e^{-ikq}\psi_1(x')={e^{i\kappa|x'|}\over |x'|}\sum\limits_{j=1}^{\infty}{f_j''(\theta)\over |x'|^{j-1}}, \ \   x'\in \R^3\setminus B_{r'},\ \  f_j''=e^{-ikq}f_j', \ \ \theta={x'\over |x'|},\eqno(4.6)$$

where the series converges absolutely and uniformly, the coefficients $f_j'$ are mentioned in (4.3).

In view of (4.4)-(4.6), similarly to the case of  Subsection 4.1,  it is sufficient to prove that  $|\psi|^2$ on $L$ uniquely determines  $f_j''(\theta)$  in (4.6) for all $j$.
This determination of $f_j''$ is completely similar to the determination of $f_j$ in Subsection 4.1.

This completes the proof of Theorem 2.1.

\vskip 2 mm
{\bf Appendix}

In this  Appendix, we  give a very simple proof that $|\psi|^2$ on $L$ uniquely determines  $f_j(\theta)$  in (3.2) for all $j$,
under the assumption that $\theta \neq k/|k|$,  where  $\psi=\psi_0+\psi_1$,  $\psi_0$  and $\psi_1$ are  solutions of (1.1) as in (1.2) and (1.3).

The determination of $f_1$ follows from (3.5).

Suppose that $f_1$,...,$f_n$ are determined, then the determination of $f_{n+1}$ is as follows.

Let

$$\psi_{1,n}(x)={e^{i\kappa|x|}\over |x|}\sum\limits_{j=1}^{n}{f_j(\theta)\over |x|^{j-1}} \ \ {\rm for}\ \  \theta={x\over |x|},   \eqno(A.1) $$

$$a_n(x,k)=|x|(|e^{ikx}+\psi_{1,n}(x)|^2-1),  \eqno(A.2)$$
$$b_n(x,k)=|x|^n(a(x,k)-a_n(x,k)),  \eqno(A.3)$$
where $x$ is as in (3.2). 

We have that

$$\eqalign{
&a(x,k)=|x|\bigl(( e^{ikx}+\psi_{1,n}(x)+{e^{i\kappa|x|}\over |x|}{f_{n+1}(\theta)\over |x|^n}+ O(|x|^{-n-2}))\cr
&\times( e^{-ikx}+\overline{\psi_{1,n}(x)}+{e^{-i\kappa|x|}\over |x|}{\overline{f_{n+1}(\theta)}\over |x|^n}+ O(|x|^{-n-2}))-1\bigr)  \cr   
& = a_n(x,k)+  {e^{i(\kappa|x|-kx)}\over |x|^n}{f_{n+1}(\theta)} + {e^{-i(\kappa|x|-kx)}\over |x|^n}\overline{f_{n+1}(\theta)}+O(|x|^{-n-1}), \cr}  \eqno(A.4) $$

$$b_n(x,k)=e^{i(\kappa|x|-kx)}{f_{n+1}(\theta)} + e^{-i(\kappa|x|-kx)}\overline{f_{n+1}(\theta)}+O(|x|^{-1}),     \eqno(A.5)$$

as $|x|\to +\infty$, uniformly in $\theta=x/|x|$, where $a$ is defined by (3.3).

Due to (A.5) and Remark 3.1,  we get
$$\eqalign{
&f_{n+1}(\theta)={1\over D}\bigl(e^{i(ky-\kappa|y|)}b_n(x,k)-e^{i(kx-\kappa|x|)}b_n(y,k)+O(|x|^{-1}) \bigr)\ \ as \ |x|\to +\infty,\cr
&D=2isin(\tau(k\theta-\kappa)),\ \  x,y\in L,\ \  y=x+\tau\theta,\ \  \theta\in\S^2,\ \  \tau>0,     \cr}  \eqno(A.6) $$  
assuming that  $D\ne 0$ for fixed $\theta$ and  $\tau$
(where  the parameter $\tau$ can be always fixed in such a way that  $D\ne 0$, under our assumption that $\theta \neq k/|k|$).

Formulas (A.1)-(A.3) and (A.6) determine  $f_{n+1}$, give the step of induction for finding all $f_j$, and complete the proof.

\vskip 2 mm
{\bf References}

\item{[  A]} F.V. Atkinson, On Sommerfeld's "Radiation Condition", Philos. Mag., Vol. XI, 645-651 (1949)

\item{[  G]} D. Gabor, A new microscopic principle,  Nature  {\bf 161}(4098), 777-778 (1948)

\item{[ JL]} P. Jonas, A.K. Louis, Phase contrast tomography using holographic measurements, Inverse Problems
{\bf 20}(1), 75-102 (2004)
\item{[ HN]} T. Hohage, R.G. Novikov, Inverse wave propagation problems without phase information,
Inverse Problems {\bf 35}(7), 070301 (4pp.)(2019)

\item{[ K]} S.N. Karp, A convergent  farfield  expansion for two-dimensional radiation functions, Communications on Pure
and Applied Mathematics {\bf 19}, 427-434 (1961)

\item{[ Kl1]} M.V. Klibanov, Phaseless inverse scattering problems in three dimensions,
SIAM J.Appl. Math. 74(2),  392-410 (2014)
\item{[ Kl2]} M.V.  Klibanov, N.A. Koshev, D.-L. Nguyen, L.H. Nguyen, A. Brettin, V.N. Astratov, A numerical method to solve a phaseless coefficient inverse problem from a single measurement of experimental data. 
SIAM J. Imaging Sci. {\bf 11}(4), 2339-2367 (2018)
\item{[ KR]} M.V. Klibanov, V.G. Romanov, Reconstruction procedures for two inverse scattering problems without the phase information, SIAM J. Appl. Math. {\bf 76}(1), 178-196 (2016)

\item{[  M]} S. Maretzke,  A uniqueness result for propagation-based phase contrast imaging from a single measurement,  
Inverse Problems {\bf 31}, 065003 (2015)
\item{[ MH]}S. Maretzke, T. Hohage,  Stability estimates for linearized near-field phase retrieval in X-ray phase contrast imaging,  SIAM J. Appl. Math.   {\bf 77}, 384-408  (2017)

\item{[ N1]} R. G. Novikov,  Inverse scattering without phase information,
S\'eminaire Laurent Schwartz - EDP et applications (2014-2015), Exp. No16, 13p.
\item{[ N2]} R.G. Novikov, Formulas for phase recovering from phaseless scattering data at fixed frequency,
Bulletin des Sciences Math\'ematiques {\bf 139}(8), 923-936 (2015)
\item{[ N3]} R.G. Novikov, Phaseless inverse scattering in the one-dimensional case,
Eurasian Journal of Mathematical and Computer Applications {\bf 3}(1), 63-69 (2015)
\item{[ N4]} R.G Novikov, Multipoint formulas for phase recovering from phaseless scattering data, J. Geom. Anal.  {\bf 31}, 1965-1991 (2021)

\item{[NSh]} R.G. Novikov,  B.L. Sharma, Phase recovery from phaseless scattering data for discrete Schrodinger operators, arXiv:2307.06041

\item{[NS1]} R.G. Novikov, V. N. Sivkin, Error estimates for phase recovering from phaseless scattering data,
Eurasian Journal of Mathematical and Computer Applications  {\bf 8}(1), 44-61 (2020)

\item{[NS2]} R.G. Novikov, V. N. Sivkin, Fixed-distance multipoint formulas for the scattering amplitude from phaseless measurements, 
Inverse Problems {\bf 38}, 025012 (2022)

\item{[ Nu]} K. Nugent,  X-ray noninterferometric phase imaging: a unified picture, Journal of the Optical Society of America. A.  {\bf 24}(2), 536-547 (2007)

\item{[  R]} V.G. Romanov,  Inverse problems without phase information that use wave interference,
Sib. Math. J. {\bf 59}(3), 494-504 (2018)

\item{[  W]} C.H. Wilcox, A generalization of theorems of Rellich and Atkinson, 
Proceedings of the American Mathematical Society  {\bf 7}(2), 271-276 (1956)

\item{[  Wo1]} E. Wolf, Three-dimensional structure determination of semi-transparent objects from holographic data, Optics
Communications  {\bf 1}(4), 153-156 (1969)

\item{[  Wo2]} E. Wolf, Determination of the amplitude and the phase of scattered fields by holography,
Journal of the Optical Society of America  {\bf 60}(1), 18-20 (1970)

\vskip 8 mm
Roman G. Novikov, CMAP, CNRS, \'Ecole polytechnique, 

Institut Polytechnique de Paris, 91128 Palaiseau, France

\& IEPT RAS, 117997 Moscow, Russia

E-mail: novikov@cmap.polytechnique.fr

\end